\documentstyle[11pt]{article}
\newlength{\standardunitlength}
\newtheorem{theorem}{Theorem} 
\newenvironment{proof}{\noindent {\sc Proof:}}{$\Box$ \vspace{2 ex}}
\begin{document}

\begin{center}
Counting Semisimple Orbits of Finite Lie Algebras by Genus
\end{center}

\begin{center}
By Jason Fulman
\end{center}

\begin{center}
Dartmouth College
\end{center}

\begin{center}
Department of Mathematics
\end{center}

\newpage

\begin{center}
Running Head: Counting Semisimple Orbits
\end{center}

\begin{center}
Mailing Address:
\end{center}

\begin{center}
Jason Fulman
\end{center}

\begin{center}
Dartmouth College
\end{center}

\begin{center}
Department of Mathematics
\end{center}

\begin{center}
6188 Bradley Hall
\end{center}

\begin{center}
Hanover, NH 03755
\end{center}

\begin{center}
email:fulman@dartmouth.edu
\end{center}

\newpage

\begin{abstract}
	The adjoint action of a finite group of Lie type on its Lie
algebra is studied. A simple formula is conjectured for the number of
semisimple orbits of a given split genus. This conjecture is proved
for type $A$, and partial results are obtained for other types. For
type $A$ a probabilistic interpretation is given in terms of card
shuffling.
\end{abstract}

\newpage

\section{Introduction}

	Let $G$ be a reductive, connected, simply connected group of Lie
type defined over an algebraically closed field of characteristic $p$. Let
$F$ denote a Frobenius map and $G^F$ the corresponding finite group of Lie
type. Suppose also that $G^F$ is $F$-split. Two semisimple elements $x,y
\in G^F$ are said to be of the same genus if their centralizers
$C_G(x),C_G(y)$ are conjugate by an element of $G^F$. It is well known in
the theory of finite groups of Lie type that character values on semisimple
conjugacy classes of the same genus behave in a unified way.

	Deriziotis \cite{Deriziotis} showed that a genus of semisimple
elements of $G^F$ corresponds to a pair $(J,[w])$, where $\emptyset
\subseteq J \subset \tilde{\Delta}, J \neq \tilde{\Delta}$ is a proper
subset of the vertex set $\tilde{\Delta}$ of the extended Dynkin
diagram up to equivalence under the action of $W$, and $[w]$ is a
conjugacy class representative of the normalizer quotient $N_W(W_J)/W_J$.
A centralizer corresponding to the data $(J,[w])$ can be obtained by twisting by $w$
the group generated by a maximal torus $T$ and the root groups $U_{\pm
\alpha}$ for $\alpha \in J$.

	Many authors (\cite{Chang}, \cite{FJ55}, \cite{FJ56}, \cite{FJ57},
\cite{Shoji}) have considered the problem of counting semisimple conjugacy
classes of $G^F$ according to genus. As emerges from their work, the number
of semisimple classes belonging to the genus $(J,[w])$ is equal to
$f(J,[w])/ |C_{N_W(W_J)/W_J} (w)|$ where $f(J,[w])$ is the number of $t$ in
a maximal torus $T$ of $G$ such that $w \cdot F(t)=t$ and the subgroup of
$W$ fixing $t$ is $W_J$. Determining $f(J,[w])$ explicitly is an elaborate
computation involving Moebius inversion on a collection of closed
subsystems of the root system.

	Let $g$ be the Lie algebra of $G$. Much less seems to be known
about the semisimple orbits of the adjoint action of $G^F$ on
$g^F$. Letting $r$ denote the rank of $G$, it is known from
\cite{Lehrer} that the number of such orbits is equal to $q^r$. By a
result of Steinberg \cite{Steinberg}, the number of semisimple
conjugacy classes of $G^F$ is also equal to $q^r$, though no
correspondence between these sets is known.

	Two semisimple elements $x,y \in g^F$ are said to be in the
same genus if $C_G(x),C_G(y)$ are conjugate by an element of
$G^F$. Experimentation with small examples such as $SL(3,5)$ suggests
that there is not an obvious relation between this decomposition of
semisimple orbits according to genus and the decomposition of
semisimple conjugacy classes of $G^F$ according to genus.

	To parametrize the genera, a semisimple element $x \in g^F$ is
said to be in the genus $(J,[w])$ where $\emptyset \subseteq J \subset
\tilde{\Delta}, J \neq \tilde{\Delta}$ if $C_G(x)$ is conjugate to the
group obtained by twisting by $w$ the group generated by a maximal
torus $T$ and the root groups $U_{\pm \alpha}$ for $\alpha \in J$.

	Lehrer \cite{Lehrer} obtained some results concerning this
parametrization. In the case where $p$ is a prime which is good and
regular (these notions are defined in Section \ref{results}) he
obtained formulae for the total number of split orbits
(i.e. $[w]=[id]$) and the total number of regular orbits
(i.e. $J=\emptyset$).

	The main conjecture of this paper is a formula for the number
of orbits in the genus $(J,[id])$ for any $J$. This formula has a
different flavor from Lehrer's formulae and counts solutions to
equations which arose in a geometric setting in Sommers' work on
representations of the affine Weyl group on sets of affine flags
\cite{Sommers}. Section \ref{results} states our conjecture, proves it
for special cases such as type $A$, and shows that it is consistent
with Lehrer's count of split orbits. Section \ref{prob} gives a
probabilistic interpretation for type $A$ involving the theory of card
shuffling. This connection is likely not as ad-hoc as it seems, given
the companion papers \cite{Fulman},\cite{Fulman2} defining card shuffling for all
finite Coxeter groups and relating it to the semisimple orbits of $G^F$ on
$g^F$.

\section{Main Results} \label{results}

	To state the main conjecture of this paper, some further
notation is necessary. Let $\Phi$ be an irreducible root system of
rank $r$ which spans the inner product space $V$. The coroots
$\check{\Phi}$ are the elements of $V$ defined as $2\alpha /
<\alpha,\alpha>$ where $\alpha \in \Phi$. Let $L$ be the lattice in
$V$ generated by $\check{\Phi}$ and set \[ \hat{L} = \{v \in V |
<v,\alpha> \in Z \ for \ all \ \alpha \in \Phi \}. \] Let $f = [\hat{L}:L]$ be the index $L$ in $\hat{L}$.
Let $\Pi =
\{\alpha_i\} \subset \Phi^+$ be a set of simple roots contained in a set of
positive roots and let $\theta$ be the highest root in $\Phi^+$. For
convenience set $\alpha_0 = - \theta$. Let $\tilde{\Pi}= \Pi \cup \{
\alpha_0 \}$. Define coefficients $c_{\alpha}$ of $\theta$ with respect to
$\tilde{\Pi}$ by the equations $\sum_{\alpha \in \tilde{\Pi}} c_{\alpha}
\alpha = 0$ and $c_{\alpha_0}=1$.

	As is standard in the theory of finite groups of Lie type,
define a prime $p$ to be bad if it divides the coefficient of some
root $\alpha$ when expressed as a combination of simple
roots. Following \cite{Lehrer}, define a prime $p$ to be
regular if the lattice of hyperplane intersections corresponding to
$\Phi$ remains the same upon reduction mod $p$. For example in type
$A$ a prime $p$ dividing $n$ is not regular.

	For subsets $S_1,S_2$ of $\tilde{\Pi}$, we write $S_1 \sim
S_2$ if there is an element $w \in W$ such that $w(S_1)=S_2$. For $S
\subset \tilde{\Pi}, S \neq \tilde{\Pi}$ and $J \subset \Pi$, the
notation $W_S \sim W_J$ means that $W_S$ is conjugate to the parabolic
subgroup $W_J$. (Although it is not always the case that $S \subset
\tilde{\Pi}, S \neq \tilde{\Pi}$ is equivalent to some $J \subset
\Pi$, it is true that $W_S$ is conjugate to some $W_J$ with $J \subset
\Pi$).

	For $\emptyset \subseteq S \subset \tilde{\Pi}, S \neq
\tilde{\Pi}$, define as in Sommers \cite{Sommers} $p(S,t)$ to be the
number of solutions ${\bf y}$ in strictly positive integers to the
equation

\[ \sum_{\alpha \in \tilde{\Pi}-S} c_{\alpha} y_{\alpha} = t. \] With these preliminaries in hand, the
main conjecture of this paper can be stated.

{\bf Conjecture 1:} Let $G$ be a reductive, connected, simply
connected group of Lie type which is $F$-split where $F$ denotes a
Frobenius automorphism of $G$. Suppose that the corresponding prime
$p$ is good and regular. Then the number of semisimple orbits of $G^F$
on ${\em g}^F$ of genus $(J,[id])$ is equal to

\[ \frac{\sum_{S \sim J} p(S,q)}{f}. \]

Remark: Sommers \cite{Sommers} studies the quantity $\sum_{S \sim J}
p(S,t)$. He shows that it can be reexpressed in either of the following
two ways, both of which will be of use to this paper.

\begin{enumerate}

\item Let $\hat{U_t}$ be the permutation representation of $W$ on the set
$\hat{L}/tL$. Let $P_1,\cdots,P_m$ be representatives of the conjugacy
classes of parabolic subgroups of $W$. Then

	\[ \hat{U_t} = \oplus_{i=1}^m [\sum_{S: W_S \sim P_i} p(S,t)]
Ind_{P_i}^W(1). \]

\item Let $A$ be a set of hyperplanes in $V=R^n$ such that $\cap_{H \in A}
H=0$. Let $L=L(A)$ be the set of intersections of these hyperplanes, where
we consider $V \in L$. Partially order $L$ by reverse inclusion and define
a Moebius function $\mu$ on $L$ by: $\mu(X,X)=1$ and $\sum_{X \leq Z \leq Y}
\mu(Z,Y)=0$ if $X<Y$ and $\mu(X,Y)=0$ otherwise. The characteristic
polynomial of $L$ is then

	\[ \chi(L,t) = \sum_{X \in L} \mu(V,X) t^{dim(X)}. \]

	For $P$ a parabolic subgroup of $W$, let $X \in L(A)$ be the fixed
point set of $P$ on $V$. Define the lattice $L^X$ to be the sublattice of
$L$ whose elements are  $\{X \cap H | H \in A \ and \ (A-H) \cap X \neq \emptyset
\}$. Let $N_W(P)$ be the normalizer in $W$ of $P$. Then

	\[ \sum_{S: W_S \sim W_J} p(S,t) = \frac{f}{[N_W(W_J):W_J]} \chi(L^X,t). \]

\end{enumerate}

	The first piece of evidence for Conjecture 1 is Theorem
\ref{consistentLehrer}, which shows consistency with Lehrer's result
\cite{Lehrer} that under the hypotheses of Conjecture 1, the total
number of split, semisimple orbits of $G^F$ on $\em{g}^F$ is equal to
$\frac{\prod_i (q+m_i)}{|W|}$, where the $m_i$ are the exponents
of $W$.

\begin{theorem} \label{consistentLehrer}

\[ \frac{1}{f} \sum_{S \subset \tilde{\Pi} \atop S \neq
\tilde{\Pi}} p(S,q) =
\frac{\prod_i (q+m_i)}{|W|} \]

\end{theorem}

\begin{proof} The left hand side is equal to $\frac{1}{f}$ times the number of
solutions in non-negative integers of the equation

\[ \sum_{\alpha \in \tilde{\Pi}} c_{\alpha} y_{\alpha} = q. \] Sommers \cite{Sommers} shows that each such
solution corresponds to an orbit of
$W$ on $\hat{L}/qL$. By Proposition 3.9 of \cite{Sommers}, the number of
fixed points of $w$ on $\hat{L}/qL$ is equal to $fq^{dim(fix(w))}$, where
$dim(fix(w))$ is the dimension of the fixed space of $w$ in its natural action on
$V$. Burnside's Lemma states that the number of orbits of a finite group $G$ on
a finite set $S$ is equal to

	\[ \frac{1}{|G|} \sum_{g \in G} Fix(g), \] where $Fix(g)$ is the number of fixed points of $g$ on $S$. Thus

\[\frac{1}{f} \sum_{S \subset \tilde{\Pi} \atop S \neq \tilde{\Pi}}
p(S,q) =  \frac{1}{|W|} \sum_{w \in W} q^{dim(fix(w))} =
\frac{\prod_i (q+m_i)}{|W|}. \] The second equality is a theorem of Shephard and Todd \cite{Shephard}.
\end{proof}

	A second piece of evidence in support of Conjecture 1 is its
truth for regular split semisimple orbits (i.e. genus
$(\emptyset,[id])$) for $G$ of classical type.

\begin{theorem} Conjecture 1 predicts that for $q=p^a$ with $p$ regular and good, the
number of regular split semisimple orbits of $G^F$ on $\em{g}^F$ is equal
to

\[ \frac{\prod_i q-m_i}{|W|}. \] This checks for types $A,B,C$ and $D$.

\end{theorem}

\begin{proof}
	Let $P_i$ be a parabolic subgroup of $W$ and $sgn$ be the
alternating character of $W$. Note by Frobenius reciprocity that
$<Ind_{P_i}^W(1),sgn>_W=0$ unless $P_i$ is the trivial subgroup, in
which case $<Ind_{1}^W(1),sgn>_W=1$. Therefore, recalling the first
formula in the remark after Conjecture 1,

\[<\hat{U}_q,sgn>_W  =  \sum_{i=1}^m <[\sum_{S:W_S \sim P_i} p(S,q)]
Ind_{P_i}^W(1),sgn>_W =  p(\emptyset,q).\] However $<\hat{U}_q,sgn>_W$ can be computed directly from its
definition. From \cite{Sommers}, if the characteristic is good then
the number of fixed points of $w$ on $\hat{L}/qL$ is equal to
$fq^{dim(fix(w))}$, where $dim(fix(w))$ is the dimension of the fixed
space of $w$ in its natural action on $V$. Thus,

\begin{eqnarray*} <\hat{U_q},sgn>_W & = & \frac{f}{|W|} \sum_{w \in W}
(-1)^{sgn(w)} q^{dim(fix(w))}\\ & = & \frac{f (-1)^r}{|W|} \sum_{w \in
W} (-q)^{dim(fix(w))}\\ & = & \frac{f \prod_i (q-m_i)}{|W|}
\end{eqnarray*} where the final equality is a theorem from
\cite{Shephard}. Comparing these expressions for $<\hat{U_q},sgn>_W$
shows that Conjecture 1 predicts that the number of regular split
semisimple orbits of $G^F$ on $\em{g}^F$ is equal to $\frac{\prod_i
q-m_i}{|W|}$.

	Let us now check this for the classical types. For type $A$, the
split semisimple orbits correspond to monic degree $n$ polynomials
which factor into distinct linear factors and have vanishing coefficient of
$x^{n-1}$. Since $p$ does not divide $n$ ($p$ is regular), by the
argument in Theorem \ref{Aworks} this is $\frac{1}{q}$ times the number of
monic degree $n$ polynomials which factor into distinct linear factors,
with no constraint on the coefficient of $x^{n-1}$. As elementary counting
shows the number of such polynomials to be $\frac{q(q-1) \cdots
(q-n+1)}{n!}$, the result follows.

	For types $B_n$ and $C_n$, split semisimple orbits correspond
to orbits of the hyperoctahedral group of size $2^n n!$ on the maximal
toral subalgebra $diag(x_1,\cdots,x_n,-x_1,\cdots,-x_n)$, where $x_i
\in F_q$ and an element $w$ of the hyperoctahedral group acts by
permuting the $x_i$, possibly with sign changes. The regular orbits
are simply those not stabilized by any non-identity $w$. To count
these orbits of the hyperoctahedral group, note first that the
hypotheses of Conjecture 1 imply that the characteristic is odd, since
$2$ is a bad prime for types $B$ and $C$. In odd characteristic the
only element of $F_q$ equal to its negative is $0$. Thus $x_1$ may be
any of the $q-1$ non-0 elements of $F_q$, $x_2$ may be any of the
$q-3$ elements of $F_q$ such that $x_2 \neq 0,\pm x_1$, and so on. As
each such hyperoctahedral orbit has size $|B_n|=|C_n|$, and the
exponents for types $B_n,C_n$ are $1,3,\cdots,2n-1$, Conjecture 1
checks for these cases.

	For type $D_n$, split semisimple orbits correspond to orbits
of $D_n$ on the maximal toral subalgebra
$diag(x_1,\cdots,x_n,-x_1,\cdots,-x_n)$ where $x_i \in F_q$ and an
element $w$ of $D_n$ acts by permuting the $x_i$, possibly with an
even number of sign changes. The regular orbits are simply those not
stabilized by any non-identity $w$. Here also one may assume odd
characteristic, as $2$ is a bad prime for type $D$. Let us consider
the possible values of $x_1,\cdots,x_n$. The first possibility is that
all $x_i \neq 0$. This can happen in $(q-1)(q-3) \cdots (q-(2n-1))$
ways, as $x_2 \neq \pm x_1$, $x_3 \neq \pm x_1,x_2$, and so on. The
second possibility is that exactly one $x_i$ is equal to $0$. As this
$i$ can be chosen in $n$ ways, the second possibility can arise in a
total of $n(q-1)(q-3) \cdots (q-(2n-3))$ ways. Thus the total number
of possible values of $x_1,\cdots,x_n$ is equal to $(q-1)(q-3) \cdots
(q-(2n-3))(q-(n-1))$. As each such orbit of $D_n$ has size $|D_n|$ and
the exponents for $D_n$ are $1,3,\cdots,2n-3,n-1$, the result follows
for type $D_n$.  \end{proof}

	As a third piece of evidence for Conjecture 1, we prove it for
$W$ of type $A$ (i.e. $SL(n,q)$). Let us make some preliminary remarks
about this case. All $p$ are good for type $A$, and it is easy to see
that if $p$ divides $n$ then $p$ is not regular for $SL(n,q)$. The
split semisimple orbits of $SL(n,q)$ on $sl(n,q)$ correspond to monic
degree $n$ polynomials $f(x)$ which factor into linear polynomials and
have vanishing coefficient of $x^{n-1}$. The genera are parameterized
by partitions $\lambda = (i^{r_i})$, where $r_i$ is the number of
irreducible factors of $f(x)$ which occur with multiplicity $i$.

\begin{theorem} \label{Aworks} Conjecture 1 holds for type
$A$. Furthermore, the number of split, semisimple orbits of $SL(n,q)$ on
$sl(n,q)$ of genus $\lambda = (i^{r_i})$ is equal to

\[ \frac{(q-1)\cdots (q+1-\sum r_i)}{r_1! \cdots r_n!}. \]

\end{theorem}

\begin{proof}
	Note that because $p$ does not divide $n$, for any $c_1,c_2$
there is a bijection between the set of split, monic polynomials with
coefficient of $x^{n-1}$ equal to $c_1$ and factorization $\lambda$,
and the set of split, monic polynomials with coefficient of $x^{n-1}$
equal to $c_2$ and factorization $\lambda$. This bijection is given by
sending $x \rightarrow x+a$ for suitable $a$. An easy combinatorial
argument shows that the number of split, monic degree $n$ polynomials
(with no restriction on the coefficient of $x^{n-1}$) of factorization
$\lambda$ is equal to

\[ \frac{q(q-1) \cdots (q+1-\sum r_i)}{r_1! \cdots r_n!}. \] Dividing
by $q$ establishes the desired count.

	To show that Conjecture 1 holds for type $A$, take $J$ of type
$\lambda$ (i.e. $W \simeq \prod S_i^{r_i}$) in the second formula in
the remark after Conjecture 1. One obtains that
	
\[ \frac{\sum_{S \sim J} p(S,t)}{f}  = 
\frac{\chi(L^X,q)}{[N_W(W_J):W_J]} =  \frac{(q-1) \cdots (q- \sum
r_i +1)}{r_1! \cdots r_n!}  \] where the formula for
$\chi(L^X,q)$ used in the second equality is Proposition 2.1 of
\cite{OrlikSolomon}.\end{proof}

\section{A Connection with Card Shuffling} \label{prob}

	The purpose of this section is to give a probabilistic proof
using card shuffling of the following identity of Lehrer \cite{Lehrer}
(which also follows from Theorem \ref{Aworks}):

\[ \sum_{\lambda=(i^{r_i}) \vdash n} \frac{q(q-1) \cdots (q-\sum r_i
+1)}{r_1!  \cdots r_n!} = \frac{q(q+1) \cdots (q+n-1)}{n!}. \] Persi
Diaconis suggested that a probabilistic interpretation might exist.

	Before doing so, we indicate the importance of card shuffling
in Lie theory and give some necessary background. For any finite
Coxeter group $W$ and $x \neq 0$, the author \cite{Fulman} defined
signed probability measures $H_{W,x}$ on $W$ as follows. For $w \in W$,
let $D(w)$ be the set of simple positive roots mapped to negative roots by $w$
(also called the descent set of $w$). Let $\lambda$ be an equivalence class of subsets
of $\Pi$ under $W$-conjugacy and let $\lambda(K)$ be the equivalence
class of $K$. Then define

\[ H_{W,x}(w) =
\sum_{K
\subseteq
\Pi-Des(w)}
\frac{|W_K|
\chi(L^{Fix(W_K)},x)}{x^n |N_W(W_K)||\lambda(K)|}.\]

	The measure $H_{S_n,x}$ arises from the theory of card
shuffling, as will be described below. It is also (expressed differently) 
related to the Poincare-Birkhoff-Witt theorem and splittings of Hochschild and cyclic
homology \cite{Loday}. There is an alternate definition of the
measures $H_{W,x}$ using the theory of hyperplane arrangements
\cite{Fulman2}. This definition leads to a concept of riffle shuffling
for any real hyperplane arrangement or oriented matroid. The measures
$H_{W,x}$ have interesting properties \cite{Fulman2}, the most remarkable of which
are:

\begin{enumerate}

\item (Convolution) $\left(\sum_{w \in W} H_{W,x}(w) w \right)
\left(\sum_{w \in W} H_{W,y}(w) w \right) = \sum_{w \in W} H_{W,xy}(w) w.
$

\item (Non-negativity) $H_{W,p}(w) \geq 0$ for all $w \in W$, provided
that $W$ is crystallographic and $p$ is a good prime for $W$.

\end{enumerate}

	There is one more (conjectural) property of $H_{W,x}$ which
should be mentioned in the context of this paper. Lehrer \cite{Lehrer}
defined a map from semisimple orbits of $G^F$ on $\em{g}^F$ to
conjugacy classes of $W$ as follows. Let $\alpha \in \em{g}^F$ be
semisimple. $G'$ simply connected implies that $C_G(\alpha)$ is
connected. Take $T$ to be an $F$-stable maximal torus in $C_G(\alpha)$
such that $T^F$ is maximally split. All such $T$ are conjugate in
$G^F$. As there is a bijection between $G^F$ conjugacy classes of
$F$-stable maximal tori in $G$ and conjugacy classes of the Weyl group
$W$, one can associate a conjugacy class of $W$ to a semisimple orbit
of $G^F$ on $\em{g}^F$. It is conjectured in \cite{Fulman} that for
$p$ good and regular, if one of the $q^r$ semisimple orbits of $G^F$
on $\em{g}^F$ is chosen uniformly at random, then the probability that
the associated conjugacy class in $W$ is a given conjugacy class $C$
is equal to the chance that an element of $W$ chosen according to the
measure $H_{W,q}$ belongs to $C$.

	The measure $H_{S_n,x}$ has an explicit ``physical'' description
when $x$ is a positive integer. This is called inverse $x$-shuffling
by Bayer and Diaconis \cite{Bayer}. Start with a deck of $n$ cards
held face down. Cards are turned face up and dealt into one of $q$
piles uniformly and independently.  Then, after all cards have been
dealt, the piles are assembled from left to right and the deck of
cards is turned face down. The chance that an inverse $q$-shuffle
leads to the permutation $\pi^{-1}$ is equal to the mass the measure
$H_{S_n,x}$ places on $\pi$.

\begin{theorem} \label{Lehrerprob}

\[ \sum_{\lambda=(i^{r_i}) \vdash n} \frac{q(q-1) \cdots (q-\sum r_i +1)}{r_1!
\cdots r_n!} = \frac{q(q+1) \cdots (q+n-1)}{n!} \]

\end{theorem}

\begin{proof}
	The right hand side is equal to the number of ways that an
inverse $q$-shuffle results in the identity. To see this, note that an
inverse $q$-shuffle gives the identity if and only if for all $r<s$,
all cards in pile $r$ have lower numbers than all cards in pile
$s$. Thus, letting $x_j$ be the number of cards which end up in the
$j$th pile, inverse $q$-shuffles resulting in the identity are in
$1-1$ correspondence with solutions of the equations $x_1 + \cdots +
x_q =n, x_j \geq 0$. Elementary combinatorics shows there to be
${q+n-1 \choose q-1}$ solutions.

	The left hand side also counts the number of ways in which an
inverse $q$-shuffle can yield the identity. As before let $x_j$ be the
number of cards which end up in the $j$th pile. The term corresponding
to $\lambda=(i^{r_i})$ counts the number of solutions to the equation
$x_1 + \cdots + x_q = n, x_j \geq 0$ and exactly $r_i$ of the $x$'s
equal to $i$. This is because such solutions are counted by the
multinomial coefficient ${q \choose q-\sum r_i, r_1, \cdots, r_n}$.
\end{proof}

	A natural problem suggested by Theorem \ref{Lehrerprob} is to
find a probabilistic interpretation of the concept of genus.

\end{document}